\documentclass[12pt]{article}
\usepackage{amsmath}
\usepackage{amssymb}
\usepackage{amsthm}
\usepackage{amsfonts}
\usepackage{amscd}
\usepackage[mathscr]{eucal}
\newcommand{\Z} {{\mathbb  Z}}
\newcommand{\Q}{{\mathbb  Q}}

\newcommand{\C}{{\mathbb  C}}

\begin{document}
\parindent  25pt
\baselineskip  10mm
\textwidth  15cm    \textheight  23cm \evensidemargin -0.06cm
\oddsidemargin -0.01cm
\title{ { A method to determine algebraically integral Cayley digraphs on finite Abelian group
 }}
\author{\mbox{}
{LI Fei}
\thanks{ \quad  Partially supported by Foundation of NSSFC(No.13CTJ006).  E-mail: cczxlf@163.com } \\
(School of Statistics and Applied Mathematics,  \\  Anhui University of Finance and Economics,\\
 Bengbu, 233030, Anhui,  P.R.China)  }

\date{}
\maketitle
\parindent  24pt
\baselineskip  10mm
\parskip  0pt

\par   \vskip 0.6cm

{\bf Abstract} \quad Researchers in the past have studied eigenvalues of Cayley digraphs or graphs.
We are interested in characterizing Cayley digraphs on a finite
Abelian group $ G $ whose eigenvalues are algebraic integers in a given
number field $K.$ And we succeed in finding a method to do so by proving Theorem 1.
Also, the number of such Cayley digraphs is computed.



\par     \vskip  1.3cm

\hspace{-0.6cm}{\bf 1 \ \ Introduction and Main Results}

\par \vskip 0.6cm

In this paper, the finite digraphs without loops and multiple edges are considered.
Let $ G $ be a group and let $ S $ be a subset of $ G $
that does not contain the identity.
Then the Cayley digraph (see [8])$ D(G, S) $ is the digraph with
vertex set $ G $ and edge set $ E(D(G, S)) = \{gh \mid hg^{-1} \in S \}. $
If $ S $ is inverse-closed, then the definition reduces
to that of a Cayley graph. When $ G $ is cyclic, $ D(G, S) $ is called circulant digraphs.
For  the  prerequisites for Cayley digraphs we refer to [5, 8].

A digraph is called algebraically integral
over a number field $ K $ if its spectrum consists of some algebraic
integers of $ K. $ Since 1974, Harary and Schwenk (see [9])
began to research integral graphs, which in fact are algebraically integral
over the rational number field $ \Q. $  There are some literature
studying integral Cayley graphs and the spectrum of Cayley digraphs(see e.g., [1][2][3][4][6][14]).

We call $ D(G, S) $ Gaussian integral if it is algebraically integral over Gaussian field $ \Q(i).$
There are several authors who have studied digraphs with Gaussian spectra(see[7][11][15]).
And in our previous results, algebraically integral circulant digraphs was discussed completely(see [13]).

In this paper,
we obtain the necessary and sufficient conditions for algebraically integral Cayley digraphs
on finite Abelian groups. Also, the number of such
algebraically integral Cayley digraphs is calculated.

Let $ n $ be the order of $ G. $ By the theorem of finite
Abelian groups, $ G $ is isomorphic to the direct sum of some cyclic groups.
In the following, we assume that $ G = C_{1} \bigoplus C_{2} \bigoplus \cdots \bigoplus C_{m}, $
where $ C_{i} = \Z/n_{i}\Z $ ($ n_{i} $ not necessary be a power of a prime number)for
some integer $ n_{i}, 1 \leq i \leq m. $ We assume without loss of generality
the fixed number field $ K $ is contained in
the $n$-th cyclotomic field $ \Q(\zeta_{n}). $
The main results are stated as follows.

{\bf A.  (see Theorem 1.) }  \ The Cayley digraph $ D(G, S) $ is algebraically integral over $ K $
if and only if $ S $ is a union of
some orbits $ Hg_{i}^{,}$s.

{\bf B.  (see Proposition 3.) }  \ For $ G = (\Z/n_{1}\Z) \bigoplus (\Z/n_{2}\Z) \bigoplus \cdots
\bigoplus (\Z/n_{m}\Z), $ there are at most $ 2^{r(G, K)} $ Cayley digraphs algebraically
integral over $ K. $

{\bf Remark. }  Based on Theorem 1, we can get a method to determine integral
and Gaussian integral Cayley digraphs on a finite Abelian group completely.

\par \vskip 0.8 cm

\hspace{-0.6cm}{\bf 2 \ \ Proof of Main Theorem }
\par  \vskip 0.2 cm

Let $ H = Gal(\Q(\zeta_{n})/K) $
be the Galois group of $ \Q(\zeta_{n})$ over $ K. $
Since $ H \subseteq Gal(\Q(\zeta_{n})/\Q) \cong (\Z/n\Z)^{\ast}, $ for each $ \sigma \in H, $
there is an element $ a \in (\Z/n\Z)^{\ast} $ such that $ \sigma(\zeta_{n}) = \zeta_{n}^{a}. $
Notice that $ n = n_{1}n_{2} \cdots n_{m}. $ For each $ n_{i}, $ we define a
group operation $ \pi_{i} $ of $ H $ on $\Z/n_{i}\Z $ by
$ \pi_{i}(\sigma, x) = \sigma(x) = ax (\text mod \ n_{i}) $
for $ x \in \Z/n_{i}\Z. $
Then we have a operation of $ H $
on $ G $ by $ \sigma (x_{1}, x_{2}, \cdots, x_{m})
= ( \pi_{1}(\sigma, x_{1}), \pi_{2}(\sigma, x_{2}), \cdots, \pi_{m}(\sigma, x_{m})) $ for
$ (x_{1}, x_{2}, \cdots, x_{m}) \in G. $ Using the orbit decomposition
formula, $ G \setminus \{0\} $ is the disjoint union of the distinct orbits,
and we can write $ G \setminus \{0\} = \bigsqcup_{i \in I}Hg_{i}, $ where $ I $ is some index
set, and the $ g_{i} $ are elements of distinct orbits.

We obtain a necessary and sufficient condition for the Cayley digraph $ D(G, S) $ to be
algebraically integral as in the following theorem.
\par \vskip 0.4 cm

 {\bf Theorem 1. } \ The Cayley digraph $ D(G, S) $ is algebraically integral over $ K $
if and only if $ S $ is a union of
some orbits $ Hg_{i}^{,}$s.

\par \vskip 0.4 cm
In the following, we give the proof of Theorem 1 in the case of $ m = 2, $ i.e., $ G $ is a direct
sum of two summands.

Let $ \widehat{G} $ denote the group of multiplicative homomorphisms from
$ G $ to $ \C^{\ast} $ (see [10]), where $ \C $ is the complex field. Then
we have isomorphism $ G \cong \widehat{G}. $
For each $ (a, b) \in G, $ there is an element $ \chi_{(a, b)} \in \widehat{G} $
such that $ \chi_{(a, b)}(s, t) =\zeta_{n_{1}}^{as} \times \zeta_{n_{2}}^{bt},
(s, t)\in G. $ The spectrum of Cayley digraph $  D(G, S) $ is given by [3],
\begin{align*}
spec(D(G, S))  & = \{ \lambda_{(0,0)}, \lambda_{(0,1)}, \cdots, \lambda_{(0,n_{2}-1)},
\lambda_{(1,0)}, \cdots, \lambda_{(n_{1}-1,0)}, \cdots, \lambda_{(n_{1}-1,n_{2}-1)} \}, \
\end{align*}
where $ \lambda_{(a,b)} = \sum_{(s, t)\in S}\chi_{(a, b)}(s, t) = \sum_{(s, t)\in S}
\zeta_{n_{1}}^{as} \cdot \zeta_{n_{2}}^{bt}. $
\par \vskip 0.4 cm

 {\bf Proposition 1. } \ If $ S $ is a union of some orbits $ Hg_{i}^{,}$s,
then $ D(G, S) $ is algebraically integral over $ K. $
\par \vskip 0.4 cm

{\bf Proof. } \ For each orbit $ Hg_{i}, $ we have $ \sigma Hg_{i} =
\{\sigma(hg_{i})\mid h\in H\} = Hg_{i}, \sigma \in H. $ So $ \sigma S = S, $ since
$ S $ is a union of some orbits $ Hg_{i}^{,}$s. Thus
\begin{align*}
\sigma(\lambda_{(a,b)})  & = \sum_{(s, t)\in S}\sigma
(\zeta_{n_{1}}^{as} \cdot \zeta_{n_{2}}^{bt})
 = \sum_{(s, t)\in S}\sigma(\zeta_{n_{1}}^{as}) \cdot
\sigma(\zeta_{n_{2}}^{bt}) \\
& = \sum_{(s, t)\in S}\zeta_{n_{1}}^{a\sigma(s)} \cdot \zeta_{n_{2}}^{b\sigma(t)} \
 =  \lambda_{(a,b)}, \
\end{align*}
for every element $ \sigma \in H. $ By the Galois theory (see [12]) for finite Galois extensions,
we get $ \lambda_{(a,b)} \in K. $ Notice that $ \lambda_{(a,b)} $ are algebraic integers, hence
$ D(G, S) $ is algebraically integral over $ K. $
\quad\quad\quad\quad $\blacksquare$

In the other hand, we have the following proposition.
\par \vskip 0.4 cm

 {\bf Proposition 2. } \ If $ D(G, S) $ is  algebraically integral over $ K, $
 then $ S $ is a union of some orbits $ Hg_{i}^{,}$s .
\par \vskip 0.4 cm

 Before the proof, we need some preparations.

Let $ \Gamma $ be a $(n-1)$-order square matrix with the row and column index in
$ G \setminus \{(0, 0)\}. $ Namely, we suppose $ \Gamma = (\gamma_{(a, b)(\alpha, \beta)}), $
where $ \gamma_{(a, b)(\alpha, \beta)} $ is the entry of the $(a, b)$-th row
and $(\alpha, \beta)$-th column, with $ (a, b), (\alpha, \beta)\in G \setminus \{(0, 0)\}. $
And $ \gamma_{(a, b)(\alpha, \beta)} = \chi_{(a, b)}(\alpha, \beta)
=\zeta_{n_{1}}^{a\alpha} \times \zeta_{n_{2}}^{b\beta}. $ For the matrix $ \Gamma, $
we have the following lemma.
\par \vskip 0.4 cm

 {\bf Lemma 1. } \ The matrix $ \Gamma $ is nonsingular.

\par \vskip 0.4 cm
{\bf Proof. } \ Let $ \chi $ be a nontrivial character of finite Abelian group $ G, $
we know that $ \sum_{g \in G}\chi(g) = 0. $ Then the sum of all entries in each row of
$ \Gamma $ is $ -1. $ Hence, to prove this lemma, it suffices to show the $ n \times n $ block
matrix
 $ \Gamma^{\prime} = \begin{pmatrix}
  1 & 1 \\
  1 & \Gamma
\end{pmatrix} $ is invertible.

Let $ F $ be the space of complex valued class function on $ G. $ All the elements
of $ \widehat{G} $ form an orthonormal basis of $ F $(see [10]). Suppose there exist $ n $
complex numbers $ k_{(a, b)}, (a, b)\in G $ such that $ \sum_{(a, b)\in G}k_{(a, b)}R_{(a, b)} =
(0, 0, \cdots, 0), $ where $ R_{(a, b)} $ is the $(a, b)$-th row of the matrix $ \Gamma^{\prime}. $
Then we have the class function $ \sum_{(a, b)\in G}k_{(a, b)}\chi_{(a, b)} = 0.$ So $ k_{(a, b)} = 0 $
for all $ (a, b) \in G, $ which shows that the row vectors of $ \Gamma^{\prime} $ are linearly
independent and so $ \Gamma^{\prime} $ is invertible. \quad\quad\quad\quad $ \blacksquare $

Let $ \tau $ be such a $(n-1)$-dimension column vector as
\begin{align*}
\tau  & = ( v_{(0,1)}, v_{(0,2)},
\cdots, v_{(0,n_{2}-1)}, v_{(1,0)}, \cdots, v_{(n_{1}-1,0)}, \cdots, v_{(n_{1}-1,n_{2}-1)})^{T} \
\end{align*}
with $ v_{(a,b)} = 1 $ for $ (a,b) \in S $ and $ 0 $ otherwise. It is easy to see that
\begin{align*}
\Gamma \tau  & = (\lambda_{(0,1)}, \lambda_{(0,2)}, \cdots, \lambda_{(0,n_{2}-1)},
\lambda_{(1,0)}, \cdots, \lambda_{(n_{1}-1,0)}, \cdots, \lambda_{(n_{1}-1,n_{2}-1)})^{T}. \
\end{align*}
Let $ \tau_{i} $ be the $(n-1)$-dimension column vector for the orbit $ Hg_{i} $ just as $ \tau $
for $ S. $ We denote $ W $ the vector
space $ \{ \omega \in K^{n-1} \mid \Gamma \omega \in K^{n-1}\} $ and $ V \subset K^{n-1} $
the the vector space spanned by the vectors $\{\tau_{i}, i \in I \}. $ We obtain the following lemma for
$ W $ and $ V. $
\par \vskip 0.4 cm

 {\bf Lemma 2. } \  $ W $ and $ V $ are the same vector space.

\par \vskip 0.4 cm
{\bf Proof. } \ By Proposition 1, $ \Gamma \tau_{i} \in K^{n-1}, i \in I. $
So $ V \subset W. $ Let
$ \omega \in W $ and  $ \omega = ( \omega_{(0,1)}, \omega_{(0,2)},
\cdots, \omega_{(0,n_{2}-1)}, \omega_{(1,0)}, \cdots, \omega_{(n_{1}-1,0)}, \cdots,
\omega_{(n_{1}-1,n_{2}-1)})^{T}, u = \Gamma \omega = ( u_{(0,1)}, u_{(0,2)},
\cdots, u_{(0,n_{2}-1)}, u_{(1,0)}, \cdots, u_{(n_{1}-1,0)}, \cdots, u_{(n_{1}-1,n_{2}-1)})^{T}. $
First, we show that $ u_{(a,b)} = u_{(c,d)} $ if
$ (a,b), (c,d) $ in the same orbit $ Hg_{i}. $  Because $ (a,b), (c,d)
\in  Hg_{i}, $ there exist an element $ \sigma \in H $ such that $ \sigma(a,b) = (c,d), $
namely, $ \sigma(a) = c (\text mod \ n_{1}) $ and $ \sigma(b) = d (\text mod \ n_{2}). $
In fact,
\begin{align*}
u_{(a,b)} = \sigma(u_{(a,b)})  & = \sigma(\sum_{(k, l)\in G\setminus\{0\}}
\omega_{(k, l)}\chi_{(a, b)}(k,l))
 = \sum_{(k, l)\in G\setminus\{0\}}\omega_{(k, l)}\sigma(\chi_{(a, b)}(k,l)) \\
& = \sum_{(k, l)\in G\setminus\{0\}}\omega_{(k, l)}\sigma(\zeta_{n_{1}}^{ak}
\cdot \zeta_{n_{2}}^{bl}) \
 = \sum_{(k, l)\in G\setminus\{0\}}\omega_{(k, l)}
\zeta_{n_{1}}^{\sigma(a)k}\cdot \zeta_{n_{2}}^{\sigma(b)l} \\
& = \sum_{(k, l)\in G\setminus\{0\}}\omega_{(k, l)}
\zeta_{n_{1}}^{ck}\cdot \zeta_{n_{2}}^{dl} = u_{(c,d)},
\end{align*}
which implies that $ \Gamma (W) \subset V. $ Notice that the matrix $ \Gamma $
is nonsingular by Lemma 1 and $ V \subset W. $ Hence $ dimW = dimV, $ and $ W = V. $
\quad\quad\quad\quad $ \blacksquare $

Now it comes to prove Proposition 2.

{\bf Proof \ of \ Proposition 2.} \quad Since $ D(G, S) $ is algebraically integral over
$ K, \Gamma\tau \in K^{n-1}. $ We have $ \tau \in W. $ By Lemma 2, $ \tau \in V $ and
$ \tau = \sum_{i \in I}c_{i}\tau_{i} $ for some coefficients $ c_{i} \in K. $ By the construction
of $ \tau $ and $ \tau_{i}^{,} $s, we conclude that $ S $ is the union of the
$ Hg_{i}^{,} $s with $ c_{i} = 1. $ The proof is completed. \quad\quad\quad\quad $ \blacksquare $
\par \vskip 0.2 cm

 {\bf Remark } Merging Proposition 1 and Proposition 2 together,
we have Theorem 1 in the case of $ m = 2. $
By [3],
$  spec(D(G, S)) = \{\lambda_{g} \mid g= (g_{1}, g_{2}, \cdots, g_{m}) \in G \}, $
where $ \lambda_{g} = \sum_{(s_{1}, s_{2}, \cdots, s_{m})\in S}
\zeta_{n_{1}}^{g_{1}s_{1}} \zeta_{n_{1}}^{g_{2}s_{2}} \cdot \zeta_{n_{2}}^{g_{m}s_{m}}. $
So far, we have given the ideas and methods to the proof of Theorem 1. It is not hard to find
these ideas and methods can be applied to the general case. In another words,
with the formula of spectrum, we can prove Theorem 1 for any $ m \in \Z^{+} $ ($ \Z^{+}, $
the set of all positive integers) in the same way as above. So we obtain the result of Theorem 1.
\par \vskip 0.4 cm

\hspace{-0.6cm}{\bf 3 \ \ Calculating the number of algebraically integral Cayley digraphs }
\par  \vskip 0.2 cm

Let $ G_{n_{i}} (1\leq i \leq m)$ be the set of all the orbits of
$ \Z/n_{i}\Z $ under $ H. $
Denote $ P_{G} $ the collection of Cartesian product $  \{ P = p_{1} \times p_{2} \times
\cdots \times p_{m} \mid \ P \neq 0, p_{i} \in G_{n_{_{i}}} \}. $ For a Cartesian product
$ P \in P_{G}, $ choose one element $ \rho \in P, \rho = (a_{1}, a_{2}, \cdots, a_{m}). $
Define $ \Q(P) $ the cyclotomic field
\begin{align*}
 \Q(\zeta_{n_{1}}^{a_{1}}, \zeta_{n_{2}}^{a_{2}},
\cdots, \zeta_{n_{m}}^{a_{m}}) = \Q(\zeta_{n_{1}}^{a_{1}}) \cdot \Q(\zeta_{n_{2}}^{a_{2}})
\cdots \Q(\zeta_{n_{m}}^{a_{m}}), \
\end{align*}
$ [\Q(P) : \Q(P)\cap K] $ the dimension of $ \Q(P) $ as vector
space over $ \Q(P)\cap K $ and $ | P | $ the cardinal number of $ P. $
It is easy to see that $ \Q(P) $ is well-defined. By Galois theory (see [12]) in number field, we have
$ | P | = \prod_{i=1}^{m}[\Q(\zeta_{n_{i}}^{a_{i}}) : \Q(\zeta_{n_{i}}^{a_{i}}) \cap K]. $

Let $ f $ be the map \
$ Gal(\Q(P)/\Q(P)\cap K) \rightarrow
  \prod_{i=1}^{m} Gal(\Q(\zeta_{n_{i}}^{a_{i}})/\Q(\zeta_{n_{i}}^{a_{i}})\cap K) $
by restriction, namely, $ \sigma \mapsto (\sigma \mid_{\Q(\zeta_{n_{1}}^{a_{1}})},
\sigma \mid_{\Q(\zeta_{n_{2}}^{a_{2}})}, \cdots, \sigma \mid_{\Q(\zeta_{n_{m}}^{a_{m}})}). $
It is easy to see that $ f $ is injective. So every orbit contained in the Cartesian product
$ P $ has $ [\Q(P) : \Q(P)\cap K] $ elements, which implies that $ P $ can be divided equally into
$ \frac{| P |}{ [\Q(P) : \Q(P)\cap K]} $ orbits under $ H. $ Totally, we get the following Lemma.
\par \vskip 0.4 cm

 {\bf Lemma 3. } \ Under the operation of $ H, $ the Cartesian product
$ P $ described as above is divided into
$ \frac{\prod_{i=1}^{m}[\Q(\zeta_{n_{i}}^{a_{i}}) : \Q(\zeta_{n_{i}}^{a_{i}}) \cap K]}{ [\Q(P) : \Q(P)\cap K]} $ orbits.
\par \vskip 0.4 cm

Let $ r(G, K) = \Sigma_{P \in P_{G}} \frac{\prod_{i=1}^{m}[\Q(\zeta_{n_{i}}^{a_{i}})
 : \Q(\zeta_{n_{i}}^{a_{i}}) \cap K]}{ [\Q(P) : \Q(P)\cap K]}. $  By Lemma 3, $ r(G, K) $ is the
orbits number of group operation of $ H $ on $ G. $ So we obtain the following proposition by
Theorem 1.
\par \vskip 0.2 cm

{\bf Proposition 3. } \ For $ G = (\Z/n_{1}\Z) \bigoplus (\Z/n_{2}\Z) \bigoplus \cdots
\bigoplus (\Z/n_{m}\Z), $ there are at most $ 2^{r(G, K)} $ Cayley digraphs algebraically
integral over $ K. $

\par \vskip 0.4 cm

\hspace{-0.8cm} {\bf References }
\begin{description}

\item[1.] RC. Alperin, BL. Peterson, Integral Sets and Cayley Graphs of Finite Groups.
Electron. J. Combin, 2012,19: \#P44.

\item[2.] J.Sander, T.Sander,  The Exact Maxamal Energy of Integral Circulant
Graphs with Prime Power Order.  Contributions to Discrete Mathematics,  2013, 8(2):19-40.

\item[3.] K. Babai, Spectra of Cayley Graphs. Journal of Combinatorial Theory,
Series B, 1979, 27: 180-189.

\item[4.] K. Balinska, D.Cvetkovic, Z.Radosavljevic, S.Simic, D.Stevanovic,
A Survey on Integral Graphs. Univ. Beograd. Publ. Elektrotehn. Fak. Ser. Mat,
2002,13: 42-65.

\item[5.] N. Biggs, Algebraic Graph theory. Amsterdam: North-Holland, 1985.

\item[6.]W. G.Bridges, R. A.Mena, Rational G-matrices with rational eigenvalues.
Journal of Combinatorial Theory, Series A, 1982,32(2): 264-280.

\item[7.] F. Esser, F. Harary, Digraphs with real and Gaussian spectra.
Discrete Appl. Math, 1980, 2: 113-124.

\item[8.] C. Godsil, G. Royle, Algebraic Graph theory. New York: Springer-Verlag, 2001.

\item[9.] F. Harary, A. Schwenk, Which graphs have integral spectra,
in: R. Bari, F. Harary (Eds.), Graphs and Combinatorics. Berlin: Springer-Verlag, 1974.

\item[10.] J. P. Serre, Linear Representations of Finite Groups.
New York: Springer-Verlag, 1977.

\item[11.] Y. Xu, J. X. Meng, Gaussian integral circulant digraphs.
Discrete Mathematics, 2011, 311: 45-50.

\item[12.] S. Lang, Algebra. 3rd edition, New York: Springer-verlag, 2002.

\item[13.] Fei Li, Circulant Digraphs Integral over Number Fields.
Discrete Mathematics, 2013, 313: 821-823.

\item[14.] W. So, Integral circulant graphs, Discrete Mathematics, 2005, 306: 153-158.

\item[15.] Xiangdong Hou, On the G-Matrices with Entries and Eigenvalues in $\Q(i)$,
Graphs and Combinatorics, 1992, 8:53-64. \ \

\end{description}

\end{document}